\documentclass[12pt]{amsart}
\usepackage{amsmath,amssymb,amsfonts,amscd,latexsym}
\newtheorem{defn0}{Definition}[section]
\newtheorem{prop0}[defn0]{Proposition}
\newtheorem{thm0}[defn0]{Theorem}
\newtheorem{lemma0}[defn0]{Lemma}
\newtheorem{corollary0}[defn0]{Corollary}
\newtheorem{example0}[defn0]{Example}
\newtheorem{remark0}[defn0]{Remark}
\newtheorem{conjecture0}[defn0]{Conjecture}

\newenvironment{proposition}{\begin{prop0}}{\end{prop0}}
\newenvironment{theorem}{\begin{thm0}}{\end{thm0}}

\newenvironment{corollary}{\begin{corollary0}}{\end{corollary0}}
\newenvironment{example}{\begin{example0}\rm}{\end{example0}}
\newenvironment{remark}{\begin{remark0}\rm}{\end{remark0}}

\newcommand{\propref}[1]{Proposition~\ref{#1}}
\newcommand{\thmref}[1]{Theorem~\ref{#1}}

\newcommand{\corref}[1]{Corollary~\ref{#1}}

\def\max{{\bf m}}                   
\def\res{{\bf k}}                   

\begin{document}
\title[Computation of Ratliff-Rush closure]{{\bf
On the  Computation of Ratliff-Rush closure}}
\author[Juan Elias ]{
Juan Elias ${}^{*}$ }
\thanks{${}^{*}$Partially supported by DGICYT BFM2001-3584, and \\
\indent Commutative Algebra program of the  MSRI, Berkeley.}
\address{Departament d'\`{A}lgebra i Geometria \newline
 \indent Facultat de Matem\`{a}tiques\newline
 \indent Universitat de Barcelona\newline
 \indent Gran Via 585, 08007
Barcelona, Spain} \email{{\tt elias@mat.ub.es}}
\date{\today}

\begin{abstract}
Let $R$ be a Cohen-Macaulay local ring with maximal ideal $\max$.
In this paper we present a procedure for computing the
Ratllif-Rush closure of  a  $\max-$primary ideal $I \subset R$.
\end{abstract}

\maketitle

\baselineskip 15pt
\medskip
\section*{Introduction}

Let $R$ be a Cohen-Macaulay Noetherian local ring of dimension $d
\ge 1$ with maximal ideal $\max$ and residue field $\res$ that we
may assume infinite.
Given a $\max$-primary  ideal $I\subset R$
in  \cite{RR78} the Ratliff-Rush closure of $I$ is defined by
$\tilde{I}=\bigcup_{k \ge 1} (I^{k+1}:I^k)),$ and it holds that
$$
\tilde{I} = \bigcup_{k \ge 1} (I^{k+1}:(x_1^k,\cdots,x_d^k))
$$
where $x_1,\cdots,x_d$ is a minimal reduction of $I$.

Although Ratliff-Rush behaves bad under most of the
basic operations of commutative algebra it is a basic tool
in the study of the Hilbert functions of primary ideals, see for
example \cite{RS02} and its reference list.

Shah defined in \cite{Sha91} a finite chain of ideals 
between $I$ and
its integral closure $\overline I$
$$
I \subset I_{[d]} \subset \cdots \subset  I_{[1]} \subset \overline{I}
$$
$ I_{[i]}$ is the $i$-coefficient ideal of $I$, and 
$ I_{[d]}=\tilde I$ the Ratliff-Rush closure of $I$.
Few results are known about the explicit computation of coefficient ideals.
Ciuperc$\breve{\text{a}}$ in \cite{Ciu01} computed
 the first coefficient ideal of an ideal $I\subset R$,
$R$ is an $(S_2)$ ring,  by considering the $S_2$-ification
of the extended Rees algebra of $I$.

\medskip
The aim of this paper is to present an algorithm for the computation of
Ratliff-Rush closure.
In the first section we prove some results on superficial sequences
that enable us to describe, in the section two, an algorithm
to compute Ratliff-Rush closure.
We end the paper with some explicit computations of   the
Ratliff-Rush closure of ideals using the algorithm of this paper.

\medskip
We will use freely \cite{BH93} as a general reference for the
algebraic concepts appearing in this paper.
The computations of this paper are performed by using CoCoA,
\cite{CocoaSystem}.

\medskip
We thank Ciuperc$\breve{\text{a}}$ for the useful comments on a previous
version of this paper.
We also  thank M.E. Rossi and  W.  Vasconcelos for pointing us that \cite[Corollary 3.4]{RTV02}
holds also for $\max-$primary ideals.

\bigskip
\section{On  superficial sequences}

Let $I$ be an $\max$-primary  ideal of $R$.
We denote by $gr_I(R)=\oplus_{k\ge 0} I^{k}/I^{k+1}$
the associated graded ring of $I$, and by $l(I)$
the analytic spread of $I$.

Let $h_{I}(n)=length_R(R/I^{n+1})$ be the Hilbert-Samuel function
of $I$, $n\in \mathbb N$. Hence there exist integers $e_j(I) \in
\mathbb Z$ such that
$$
p_{I}(X)=\sum_{j=0}^{d} (-1)^j e_j(I) \binom{X+d-j}{d-j}
$$
is the  Hilbert-Samuel polynomial of $I$, i.e. $h_{I}(n)=p_{I}(n)$
for $n \gg 0$. The integer $e_j(I)$  is the $j-$th
Hilbert coefficient of $I$, $j=0,\cdots , d$.
Shah proved that coefficient ideals are the largest ideals $I_{[t]}$ containing $I$
and such that:
\begin{enumerate}
\item[(i)] $e_i(I)=e_i(I_{[t]})$ for $i=0,\cdots,t$,
\item[(ii)] $I \subset I_{[d]} \subset \cdots \subset  I_{[1]}
\subset \overline{I}$
\end{enumerate}

\noindent
where $ \overline{I}$ is the integral closure of $I$, \cite{Sha91}.
Notice that $ \tilde I $ is the largest ideal containing $I$
and such that $e_i(I)=e_i(\tilde I)$ for $i=0,\cdots, d$.

We say that $x\in I$ is a superficial element of $I$ if there exists
an integer $k_0$ such that $(I^{k+1}:x)=I^k$ for
$k\ge k_0$.
Since the residue field is infinite it hold:
\begin{enumerate}
\item  a  set elements $x_1,\cdots, x_d \in I$,
such that their  cosets
$\overline{x}_1,\cdots, \overline{x}_d\in I/\max I$ are generic,
form a superficial sequence $x_1,\cdots,x_d$ of $I$, i.e.
$x_i$ is a superficial element of $I/(x_1,\cdots,x_{i-1})$
for $i=1,\cdots,d$,
\item if $x_1,\cdots,x_d \in I$ is a set of elements
such that
$$
\text{ length }_R\left(\frac{R}{(x_1,\cdots,x_d)}\right)=e_0(I),
$$
where $e_0(I)$ is the multiplicity of $I$, then $x_1,\cdots,x_d$ a superficial sequence
of the ideal $I$,
\item if $x_1,\cdots,x_d$ is a superficial sequence
of $I$ then $J=(x_1,\cdots,x_d)$ is a minimal reduction of $I$, \cite{Swa94}.
\end{enumerate}

\noindent
Given a superficial element $x$ of $I$ if we write $\overline{I}=I/(x)$
then it is well known that
$$
p_{\overline{I}}(X)=p_{I}(X)-p_{I}(X-1)=\sum_{j=0}^{d-1} (-1)^j e_j(I)
\binom{X+d-2-j}{d-1-j}
$$
\noindent
in particular $e_i(\overline{I})=e_i(I)$ for $i=0,\cdots,d-1$.
We define the postulation number $pn(I)$ of $I$ as the smallest
integer $n$ such that $h_I(t)=p_I(t)$ for all $t \ge n$.
Given a superficial sequence $x_1,\cdots,x_d$ of $I$
we denote by $pn(I;x_1,\cdots,x_d)$ the maximum  among
$pn(I)$ and $pn(I/(x_i))$, $i=1,\cdots,d$.

\bigskip

\begin{proposition}
\label{super}
Let $I$ be a $\max-$primary ideal of $R$ and  $x$  a superficial element of $I$.
 We denote by
$\overline{I}=I/(x)$ the ideal of $\overline{R}=R/(x)$.
For all $k \ge pn(I;x) +1$
it holds
$$
(I^{k+1}:x)=I^k.
$$
\end{proposition}
\begin{proof}
Let us consider the exact sequence
$$
0 \longrightarrow
\frac{(I^{k+1}:x)}{I^k} \longrightarrow
\frac{R}{I^k}
\overset{.x}{\longrightarrow}
\frac{R}{I^{k+1}}\longrightarrow
\frac{\overline{R}}{\overline{I}^{k+1}}\longrightarrow 0,
$$
so
$$
\text{ length }_R\left(\frac{(I^{k+1}:x)}{I^k}\right)=
h_{I}(k-1)-h_I(k)+h_{\overline{I}}(k).
$$

\noindent
If $k \ge pn(I;x)+1$ then we have
that
$h_{I}(k)=p_I(k)$,
$h_{I}(k-1)=p_I(k-1)$
and $h_{\overline{I}}(k)=p_{\overline{I}}(k)$, so
$$
\text{ length }_R\left(\frac{(I^{k+1}:x)}{I^k}\right)=
p_{I}(k-1)-p_I(k)+p_{\overline{I}}(k).
$$

On the other hand, since $x$ is a superficial element of $I$
we have that $p_{\overline{I}}(X)=p_{I}(X)-p_{I}(X-1)$ then
$(I^{k+1}:x)=I^k$ for all $k \ge pn(I;x)+1$.
\end{proof}

\bigskip
Notice that for the explicit computations of coefficient ideals  it is
enough to consider the number $pn(I;,x_1,\cdots,x_d)$, \thmref{coeff} $(i)$,
but if we look for a explicit formula of the Ratliff-Rush closure
avoiding the computation of superficial sequences we have to consider
the Castelnuovo-Mumford regularity, \thmref{coeff} $(ii)$.

\medskip

Given a standard $A_0$-algebra $A=A_0\oplus A_1 \oplus \cdots$
with $A_0$ an Artin ring,
 we denote by $\text{ reg }(A)$ the
Castelnuovo-Mumford regularity of $A$, i.e. the smallest integer
$m$ such that $H^{i}_{A_{+}}(A)_n=0$
for all $i=0,\cdots,d$ and  $n \ge m-i+1$, where
$A_{+}=A_1\oplus \cdots$ the irrelevant ideal of $A$.

We denote by $f: \mathbb N ^2 \longrightarrow \mathbb N$
the  numerical function defined by

$$
f(e,d)=
\begin{cases}
e -1 & \text {if } d=1 \\
e^{2(d-1)!-1}(e-1)^{(d-1)!} & \text {if } d\ge 2
\end{cases}
$$

\noindent
Rossi, Trung and Valla prove that $f(e,d)$ is an upper bound of
the Castelnuovo-Mumford regularity of the associated
graded ring of $I$, see \propref{RTV}.

Given a minimal reduction $J$ of $I$ we denote by $r_J(I)$
the reduction number of $I$ with respect to $I$, i.e.
the smallest integer $r$ such that $I^{r+1}=J I^r$.

In the next result we relate some of the numerical characters
that we already defined in this paper.

\bigskip
\noindent
\begin{proposition}
\label{RTV}
Let $R$ be a Cohen-Macaulay local ring of dimension $d\ge 1$.
Let $I$ be a $\max$-primary ideal of $R$ and $J=(x_1,\cdots, x_d)$ a minimal reduction of $I$.
Then
\begin{enumerate}
\item[(i)]
$r_J(I) \le  \text{ reg }(gr_{I}(R)) \le f(e_0(I), d)$, and
\item[(ii)]
$pn(I;x_1,\cdots, x_d) \le  f(e_0(I), d) +1$.
\end{enumerate}
\end{proposition}
\begin{proof}
$(i)$ The first inequality comes from \cite[Proposition 3.2]{Tru87}, see also \cite[Theorem 18.3.12]{BSLC}.
The second inequality is due to Rossi, Trung and Valla,  \cite[Corollary 3.4]{RTV02} .

\noindent
$(ii)$
Notice that from Serre's formula, \cite[Theorem 4.4.3]{BH93}, and the right hand side inequality in  $(i)$
we have that
$pn(I) \le f(e_0(I), d)+1$ and $pn(I/(x_i)) \le  f(e_0(I/(x_i)), d-1)+1$, $i=1,\cdots,d$.
Since
$e_0(I)=e_0(I/(x_i))$ and $f(e,d-1) \le f(e,d) $,
we get the claim.
\end{proof}

\bigskip
\noindent
Notice that in \cite[ Corollary 3.4]{RTV02}  the right hand side inequality in $(i)$ of the above result is proved
for the maximal ideal $I=\max$, but the proof holds also
for general $\max-$primary ideals.

\bigskip
\begin{corollary}
\label{super2}
Let $x$ be a superficial element of $I$.
For all $k \ge f(e_0(I), d) +2$ we have
$$
(I^{k+1}:x)=I^k.
$$
\end{corollary}
\begin{proof}
It is a consequence of \propref{RTV} $(ii)$ and \propref{super}.
\end{proof}

\bigskip
\section{An algorithm for computing  Ratliff-Rush closure}

\medskip
In this section we compute explicitly Ratliff-Rush closure by using
\propref{super} and \corref{super2}.
We consider  the increasing ideal chain
$$
\mathcal L_{1} \subset \mathcal L_{2} \subset \cdots \subset
\mathcal L_{k} \subset \cdots
$$
where
$$
\mathcal L_{k} = (I^{k+1}:(x_1^k,\cdots,x_d^k)).
$$
Notice that   $\tilde I = \bigcup_{k \ge 1} \mathcal L_{k}$
is the Ratliff-Rush closure of $I$.

\bigskip
\begin{theorem}
\label{coeff}
Let $R$ be a Cohen-Macaulay  local ring of dimension $d
\ge 1$.
Let $I$ be an $\max$-primary  ideal of $R$ and let
$x_1,\cdots,x_d$ be a superficial sequence of $I$.
\begin{enumerate}
\item[(i)]
For all $k \ge pn(I;x_1,\cdots,x_d)+1$ it holds that
$$
\tilde I=(I^{k+1}:(x_1^k,\cdots,x_d^k)).
$$
\item[(ii)]
For all $k\ge (d+1)(f(e_0(I))+2)$ it holds that
$$
\tilde I = (I^{k+1}:I^k).
$$
\end{enumerate}
\end{theorem}
\begin{proof}
$(i)$
We have to prove that for all $k\ge pn(I;x_1,\cdots,x_d)+1$ it holds
$\mathcal L_{k}=\mathcal L_{k+1}.$
Notice that for all $n \ge 1$ we have
 $\mathcal L_{n} \subset \mathcal L_{n+1}$ so
 we only need to prove
$\mathcal L_{k+1} \subset \mathcal L_{k}$.
Given $a  \in \mathcal L_{k+1} $ we have
$a x_i^{k+1} = x_i (a x_i^{k}) \in I^{k+2}$, for all $i=1\cdots, d$.
Since $k\ge pn(I;x_1,\cdots,x_d)+1$ from \propref{super} we get
$a x_i^{k} \in I^{k+1}$ for all $i=1, \cdots, d$, so
$a \in \mathcal L_{k}$.

\noindent$(ii)$
Notice that $J=(x_1,\cdots,x_d)$ is a  minimal reduction of $I$
so for all $k\ge r_J(I)$ 
$$
I^{(d+1)k}=I^{dk}(x_1^k,\cdots,x_d^k).
$$

\noindent
From  \propref{RTV}  we have
 that $r_J(I) \le  reg(gr_I(R)) \le f(e_0(I),d)$.
Let $n \ge f(e_0(I),d)+2$ be an integer and let
$a\in \tilde{I}$ be an element of the Rattlif-Rush closure of $I$.
Hence from  $(i)$ we have $a x^{[k]} \subset I^{k+1}$
and since $I^{(d+1)k}=I^{dk}(x_1^k,\cdots,x_d^k)$ we get
$$
a I^{(d+1)k} \subset a I^{dk}(x_1^k,\cdots,x_d^k) \subset
I^{(d+1)k+1}.
$$
In particular we have
$a \in (I^{(d+1)k+1}: I^{(d+1)k})$,
since by definition $(I^{(d+1)k+1}: I^{(d+1)k}) \subset \tilde{I}$
we get the claim.
\end{proof}

\bigskip
From the last  result we deduce that the problem of computing the Ratliff-Rush closure
can  be reduced  to the computation of  the postulation number of $I$ and its quotients
$I/(x_i)$, $i=1,\cdots,d$.
Next we  recall how to compute these numbers.

We denote by $PS_I(X) \in {\mathbb Z}[[X]]$ the Poincar{\'e}  series of $I$
$$
PS_I(X) = \sum_{i\ge 0} length_R\left(\frac{I^i}{I^{i+1}}\right) X^i
$$
\noindent
it is known that there exists a degree $s$ polynomial
$f(X)=\sum_{i=0}^{s} a_i X^  i \in \mathbb Z [X]$
such that
$$
PS_I(X) =\frac{f(X)}{(1-X)^d}.
$$
\noindent
It is easy to prove that $e_0(I)=\sum_{i=0}^{s} a_i$ and that
$pn(I)=s-d$.

\begin{remark}
It is well known  that the computation of the  Poincar{\'e} series of $I$ and its quotients
$I/(x_i)$  can  be reduced to a elimination of variables process, see for example
the library { \tt primary.lib} of CoCoa, \cite{CocoaSystem}.
\end{remark}

\bigskip
\noindent
\begin{center}
{\sc An algorithm for  computing the Ratliff-Rush closure.}
\end{center}
\medskip

\begin{enumerate}
\item[\bf{Step 1.}] Compute the Poincar{\'e} series of $I$.
Then we know the multiplicity $e_0(I)$
and the postulation number $pn(I)$ of $I$.
\item[\bf{Step 2.}] Find $d$ generic elements $x_1,\cdots,x_d$ of the $\res-$vector space $I/\max I$
such that $length_R(R/(x_1,\cdots,x_d))=e_0(I)$.
Recall that  $x_1,\cdots,x_d$ is a superficial sequence of $I$ and generates  a minimal reduction of
$I$.
\item[\bf{Step 3.}]
As in {\bf{Step 1}} compute $PS_{I/(x_i)}=f_i(X)/(1-X)^{d-1}$ for $i=1,\cdots,d$.
From this and the fact $pn(I/(x_i))=deg(f_i)-(d-1)$
we can compute  $pn(I;x_1,\cdots,x_d)$.
\item[\bf{Step 4.}]
For $k \ge pn(I;x_1,\cdots,x_d)+1$  we get
$$
\tilde I =(I^{k+1}:(x_1^k,\cdots,x_d^k)).
$$
\end{enumerate}

\begin{remark}
Notice that if $I$ is a monomial ideal then {\bf{Step 4}} can be performed without
Gr{\"o}bner basis computation.
\end{remark}

\bigskip
We will show how to compute the Ratliff-Rush closure
 in some  explicit examples of \cite{Ciu01} and
\cite{RS02}.

\begin{example} Example 1.10 of \cite{RS02}.
Let $I=(x^{10},y^5,xy^4,x^8y)$ be an ideal  of $R=\res[x,y]_{(x,y)}$.
The Poincar{\'e} series of $I$ is
$$PS_I(X)=\frac{35 + 4X + 4X^2 + 4X^3 - 2X^4}{(1-X)^2},
$$
so $e_0(I)=45$ and $pn(I)=2$.
Since the length of $R/(y^5+x^{10}+x^8y,xy^4)$ is $45=e_0(I)$ we deduce that
$y^5+x^{10}+x^8y,xy^4$ is a superficial sequence of $I$.
A CoCoA  computation shows that
$$
PS_{I/(xy^4)}(X)=\frac{35 + 6X + 2X^2 + 2X^3}{1-X},
$$
and
$$
PS_{I/(y^5+x^{10}+x^8y)}(X)=\frac{35 + 6X + 4 X^2}{1-X},
$$
so $pn(I;y^5+x^{10}+x^8y,xy^4)=2$.
Then by \thmref{coeff} $(i)$ we get
$$I \subsetneq \tilde{I}=(I^4:((y^5+x^{10}+x^8y)^3,(xy^4)^3)))
=(x^{10}, y^5, xy^4, x^7y^2, x^6y^3, x^8y).$$
\end{example}

\begin{remark}
Let $x_1=y^5+x^{10}+x^8y,x_2=xy^4$ be the minimal reduction of the ideal
$I$ of the last example.
Since  $pn(I;x_1,x_2)=2$ we have that
$\tilde{I}=(I^4:(x_1^3,x_2^3))$.
On the other hand, \thmref{coeff} $(ii)$ gives that
$\tilde{I}=(I^{k+1}:I^k)$ for all $k \ge 540$,  this  is
a hard computation.
\end{remark}

\begin{example} Example 1.4 of \cite{RS02}.
Let us consider the ideal
$$I=(y^{22},x^4y^{18},x^7y^{15},x^8y^{14},x^{11}y^{11},x^{14}y^8,x^{15}y^7,x^{18}y^4,x^{22})$$
of the local ring $R=\res[x,y]_{(x,y)}$.
A similar computation as we did in the previous example
shows that $I=\tilde{I}$ and
$$
I^2 \subsetneq
\widetilde{I^2}=I^2+(x^{24}y^{20}, x^{20}y^{24}).
$$
\end{example}

\begin{example} Example 3.3  of \cite{Ciu01}.
\label{cuip}
Let us consider the ideal
$$I=(x^8,x^3y^2,x^2y^4,y^8)$$
of the local ring  $ R=\res[x,y]_{(x,y)} $.
A similar computation as before shows that $I=\widetilde{I}$.

\noindent
Ciuperc$\breve{\text{a}}$ in \cite{Ciu01} computed
 the first coefficient ideal of $I$:
 $$
I=\widetilde{I} \subsetneq I_{[1]}=(x^8, x^3y^2, x^2y^4, xy^6,  y^8).
$$

\end{example}

\baselineskip 8pt



\begin{thebibliography}{10}

\bibitem{BSLC}
M.P. Brodman and R.Y. Sharp.
\newblock {\em Local Cohomology}, volume~60 of {\em Cambridge {S}tudies in
  {A}dvanced {M}athematics}.
\newblock Cambridge University Press, 1998.

\bibitem{BH93}
W.~Bruns and J.~Herzog.
\newblock {\em Cohen-{M}acaulay Rings}.
\newblock Cambridge studies in advanced mathematics. Cambridge University
  Press, 1993.

\bibitem{CocoaSystem}
A.~Capani, G.~Niesi, and L.~Robbiano.
\newblock {\em {CoCoA, a system for doing Computations in Commutative
  Algebra}}.
\newblock {Available via anonymous ftp from {\tt cocoa.dima.unige.it}}, {4.0}
  edition, 2000.

\bibitem{Ciu01}
C.~Ciuperc$\breve{\text{a}}$.
\newblock First coefficient ideals and the {$S_2$}-ification of a {R}ees
  algebra.
\newblock {\em J. of Alg.}, 242:782--794, 2001.

\bibitem{RR78}
L.~Ratliff and D.~Rush.
\newblock Two notes on reductions of ideals.
\newblock {\em Indiana Univ. Math. J.}, 27:929--934, 1978.

\bibitem{RS02}
M.~E. Rossi and I.~Swanson.
\newblock Notes on the behaviour of the {R}atliff-{R}ush filtration.
\newblock {\em Preprint}, 2002.

\bibitem{RTV02}
M.E. Rossi, N.V. Trung, and G.~Valla.
\newblock Castelnuovo-{M}umford regularity and extended degree.
\newblock {\em Preprint}, 2002.

\bibitem{Sha91}
K.~Shah.
\newblock Coefficient ideals.
\newblock {\em Trans. A.M.S.}, 327(1):373--384, 1991.

\bibitem{Swa94}
I.~Swanson.
\newblock {A} note on analytic spread.
\newblock {\em Comm. in Algebra}, 22:407--411, 1994.

\bibitem{Tru87}
N.V. Trung.
\newblock Reduction exponent and degree bound for the defining equations of
  graded rings.
\newblock {\em Proc. A.M.S.}, 101-2:229--236, 1987.

\end{thebibliography}

\end{document}